\documentclass[a4paper,12pt]{article}
\usepackage{times, url}
\usepackage{graphicx,amsmath,amssymb,amsthm,amsfonts,bm,float,cite}
\usepackage{kpfonts}
\usepackage[T1]{fontenc}

\newtheorem{definition}{Definition}[section]
\usepackage[none]{hyphenat}
\usepackage[center]{caption}
\usepackage[pdfencoding=auto]{hyperref}
\usepackage{bookmark}
\usepackage{mathtools}
\DeclarePairedDelimiter\ceil{\lceil}{\rceil}

\begin{document}

\begin{center}
{\Large \bf GENERAL UNCROSSING COVERING PATHS\\ \vspace{1mm}
INSIDE THE AXIS-ALIGNED BOUNDING BOX}
\vspace{12mm}

{\Large Marco Rip\`a}
\vspace{16mm}
\end{center}

\noindent {\bf Abstract:} \sloppy Given the finite set of $n_1 \cdot n_2 \cdot \ldots \cdot n_k$ points $G_{n_1,n_2,\ldots,n_k} \subset \mathbb{R}^k$ such that $n_k \geq \cdots \geq n_2 \geq n_1 \in \mathbb{Z}^+$, we introduce a new algorithm, called M$\Lambda$I, which returns an uncrossing covering path inside the minimum axis-aligned bounding box $[0,n_1-1]  \times [0,n_2-1] \times \cdots \times [0,n_k-1]$, consisting of $3 \cdot \prod_{i=1}^{k-1} n_i-2$ links of prescribed length $n_k-1$ units. Thus, for any $n_k \geq 3$, the link length of the covering path provided by our M$\Lambda$I-algorithm is smaller than the cardinality of the set $G_{n_1,n_2,\ldots,n_k}$. Furthermore, assuming $k>2$, we present an uncrossing covering path for $G_{3,3,\ldots,3}$, consisting of $20 \cdot 3^{k-3}-2$ straight-line edges that are $2$ units long each, which is constrained by the axis-aligned bounding box $\left[0,4-\sqrt{3}\right] \times \left[0,4-\sqrt{3}\right] \times [0, 2]^{k-2}$. 
\vspace{8mm}


\section{Introduction} \label{sec:Intr}

In the present paper, we will study constrained optimization problems that are strongly related to three-dimensional integrated circuits (3D ICs) design and could have applications in robot manufacturing. In order to compactly describe these fundamental problems, and extend our results to higher dimensions, we shall do well to begin with a few basic definitions.

\begin{definition} \label{def1}
Let $W$ and $Z$ be two sets, and let us denote the Cartesian product by the symbol ``$\times$''. Thus, $W \times Z \coloneqq \{(w,z): w \in W \wedge z \in Z\}$, since it represents the set of all points $(w,z)$, where $w \in W$ and $z \in Z$.
\end{definition}

\begin{definition} \label{def2}
Let $n_1,n_2,\ldots,n_k \in \mathbb{Z}^+$ be such that $n_1 \leq n_2 \leq \cdots \leq n_k$.
We define, for every positive integer $k$, $G_{n_1,n_2,\ldots,n_k} \coloneqq \{0,1,\ldots,n_1-1\} \times \{0,1,\ldots,n_2-1\} \times \cdots \times \{0,1,\ldots,n_k-1\}$ so that the grid $G_{n_1,n_2,\ldots,n_k}$ is a set of $\prod_{i=1}^{k} n_i$ points in the Euclidean space $\mathbb{R}^k$.
\end{definition}

\begin{definition} \label{def3}
For every positive integer $k$, let us denote by $\hat{B}_{n_1,n_2,\ldots,n_k} \coloneqq \{(x_1,x_2,\ldots,x_k) : x_1 \in [0,n_1-1] \wedge x_2 \in [0,n_2-1] \wedge \cdots \wedge x_k \in [0, n_k-1]\}\subset \mathbb{R}^k$ the MAABB (i.e., minimum axis-aligned bounding box), and let us define $B_{n_1,n_2,\ldots,n_k } \coloneqq \{(x_1,x_2,\ldots,x_k) : x_1 \in [0,n_1 ] \hspace{1mm} \wedge \linebreak x_2 \in [0,n_2 ] \hspace{1mm} \wedge \cdots \wedge \hspace{1mm} x_k \in [0,n_k]\} \subset \mathbb{R}^k$ as RAABB (i.e., regular axis-aligned bounding box).
\end{definition}

\begin{definition} \label{def4}
A covering path is a directed polygonal chain that visits every node of $G_{n_1,n_2,\ldots,n_k}$ exactly once, while a covering trail is a directed polygonal chain that joins every node of $G_{n_1,n_2,\ldots,n_k}$, but it can visit any of them more than once; in both cases, edges cannot be repeated and every couple of consecutive edges cannot be collinear. A covering path is uncrossing if none of its edges meets one other, and it is self-intersecting otherwise.
\end{definition}

\begin{definition} \label{def5}
A covering circuit $F_{n_1,n_2,\ldots,n_k}$ is a closed directed trail that visits every node of $G_{n_1,n_2,\ldots,n_k}$ in which the starting point is equal to the endpoint. A covering cycle for $G_{n_1,n_2,\ldots,n_k}$ is a covering circuit in which the only repeated node can be the first/last one. In particular, we call a regular covering cycle a cycle whose endpoint/starting point belongs to $G_{n_1,n_2,\ldots,n_k}$ and we call a smart covering cycle a cycle whose endpoint/starting point is a Steiner point (i.e., a point that does not belong to $G_{n_1,n_2,\ldots,n_k}$).
\end{definition}

\begin{definition} \label{def6}
The link length $h(Q)$ of a covering trail/circuit $Q_{n_1,n_2,\ldots,n_k}$ for $G_{n_1,n_2,\ldots,n_k}$ corresponds to the number of its edges, while the length of a given edge is the Euclidean distance between its two endpoints. Lastly, let us define the length classes of a covering trail/circuit as the set comprising all the lengths of its edges, denoting by $l_{1 \leq j \leq \prod_{i=1}^{k} n_i-1} (Q)$ the $j$-th element of the set mentioned above (where $l_1(Q)<l_2(Q)< \cdots <l_{j_{max}}(Q)$) and simplify the notation by omitting the subscript if $l_1(Q)=l_{j_{max}}(Q)$, so that $l(Q) \coloneqq l_1(Q)$.
\end{definition}

Thanks to the new M$\Lambda$I-algorithm, described in Section \ref{sec:2}, we solve a pivotal problem concerning uncrossing covering paths consisting of less than $\prod_{i=1}^{k} n_i-1$ links of prescribed length \cite{1}, lying entirely inside the MAABB of any $k$-dimensional finite set of $n_1 \times n_2 \times \cdots \times n_k$ points.

In Section \ref{sec:2}, we also show how, for specific cases as $n_1=n_k=3$, it is possible to shorten the link length \cite{2} of the general solution provided by the M$\Lambda$I-algorithm if we consider the bounding box $\prod_{i=1}^{k} [0, n_i]$  (i.e., the RAABB) instead of the MAABB. Moreover, referring to the grid graphs $G_{3,3}=\{0,1,2\} \times \{0,1,2\}$ and $G_{3,3,3}=\{0,1,2\} \times \{0,1,2\} \times \{0,1,2\}$, we have constructively proved in \cite{3} the existence of self-intersecting covering paths inside the MAABB (see Definition~\ref{def3}) that are formed by less than $3^k$ line segments, all belonging to the same irrational length class.

This work aims to find valid solutions for some problems that arise when we combine in different ways a few basic constraints on the paths needed to cover all the points of the set $G_{n_1,n_2,\ldots,n_k}$ (see Definition~\ref{def2}) \cite{2, 4}.

The four fundamental constraints we are interested in, are as follows:
\begin{enumerate}
\item Visit all the points of $G_{n_1,n_2,\ldots,n_k}$ with an uncrossing covering path $P_{n_1,n_2,\ldots,n_k}$ \cite{5};
\item Consider only covering paths such that all their edges belong to a unique length class $l(P_{n_1,n_2,\ldots,n_k}) \in \mathbb{R}$;
\item $P_{n_1,n_2,\ldots,n_k} \subset B_{n_1,n_2,\ldots,n_k}$;
\item The link length of $P_{n_1,n_2,\ldots,n_k}$ has to be smaller than the cardinality of the set $G_{n_1,n_2,\ldots,n_k}$. Therefore let $h(P_{n_1,n_2,\ldots,n_k})< \prod_{i=1}^{k} n_i$.
\end{enumerate}

The additional constraint that will be considered is to replace the RAABB with the MAABB \cite{3}. Thus, merely $P_{n_1,n_2,\ldots,n_k } \subseteq \hat{B}_{n_1,n_2,\ldots,n_k}$, instead of $P_{n_1,n_2,\ldots,n_k} \subset B_{n_1,n_2,\ldots,n_k}$ (as stated by the third rule above).

The goal to minimize the total (Euclidean) length of the covering path, denoted by $\lambda(P_{n_1,n_2,\ldots,n_k})$, will not be taken into account in this paper \cite{6, 7}, and it is trivial to note that it cannot be less than $\prod_{i=1}^{k} n_i-1$ units, which implies $l(P_{n_1,n_2,\ldots,n_k})=1$. Consequently, we can be interested in reducing/minimizing the link length of $P_{n_1,n_2,\ldots,n_k}$ if and only if all the previously stated conditions have been fulfilled so that $l(P_{n_1,n_2,\ldots,n_k}) \cdot h(P_{n_1,n_2,\ldots,n_k} )=\lambda(P_{n_1,n_2,\ldots,n_k})$.

Before asking if there exists any covering path which satisfies the four fundamental constraints and if $l(P_{n_1,n_2,\ldots,n_k})$ is unique or not for some $k$-tuple $(n_1,n_2,\ldots,n_k)$, let us give a valid lower bound on the link length of any covering trail for $G_{n_1,n_2,\ldots,n_k}$.

Since Reference \cite{5}, Equation 4, guarantees a lower bound for the link length of the minimal covering trail (without additional constraints as above), and considering that, by definition, it cannot be greater than the link length of the minimal covering path, $\forall k,n \in \mathbb{N}-\{0,1,2\}$ we have
\begin{equation} \label{eq1}
h(P_{n_1,n_2,\ldots,n_k} ) \geq \ceil*{3 \cdot \frac{\prod_{i=1}^{k} n_i \sum_{i=1}^{k-2} n_i+k-3}{2 \cdot n_k+n_{k-1}-3}}+k-2.
\end{equation}


\section{The M\texorpdfstring{$\bm{\Lambda}$I }- algorithm} \label{sec:2}

We introduce the general M$\Lambda$I-algorithm, which can cover $G_{n_1,n_2,\ldots,n_k}$ for any given
k-tuple $(n_1,n_2,\ldots,n_k)$ such that $k \geq 2$, under the four fundamental constraints stated in the previous section, plus the additional one. The resulting, inside the MAABB, uncrossing covering path has a total of $\prod_{i=1}^{k-1}n_i+2 \cdot \left(\left(\prod_{i=1}^{k-1} n_i \right)-1 \right)$ edges.

Hence,
\begin{equation} \label{eq2}
h(P_{n_1,n_2,\ldots,n_k} )=3 \cdot \prod_{i=1}^{k-1}n_i-2. 
\end{equation}

Since $k$ is greater than $1$ by hypothesis, the inequality
\begin{equation} \label{eq3}
3 \cdot \prod_{i=1}^{k-1}n_i-2 < \prod_{i=1}^{k}n_i 	
\end{equation}
holds for any $n_1 \leq n_2 \leq \cdots \leq n_k$ such that $n_k>2$.

Let $n_1>1$. If $n_k=2$, it is clearly possible to join the $2^k$ vertices of the grid spending $2^k-1$ links of unitary length, as explained in \cite{6}, reproducing $2^{k-2}$ times the trivial two-dimensional covering path $P_{2,2}=(0,0)$-$(0,1)$-$(1,1)$-$(1,0)$ and fixing the aforementioned $3 \cdot 2^{k-2}$ edges with $2^{k-2}-1$ additional links. Similarly, if $k=1$, then $P_{n_1}=(0,0)$-$(0,n_1) \equiv \hat{B}_{n_1}$ works.

Thus, from here on, let us assume $3 \leq n_k \wedge k \geq 2$, and the M$\Lambda$I-algorithm can always be applied, providing a valid path such that $h(P_{n_1,n_2,\ldots,n_k})+1<\prod_{i=1}^{k} n_i$.

\textbf{M$\bf{\Lambda}$I-algorithm:} Given $(n_1,n_2,\ldots,n_k) \Rightarrow G_{n_1,n_2,\ldots,n_k}$, return an uncrossing covering path $P_{n_1,n_2,\ldots,n_k} \subseteq \hat{B}_{n_1,n_2,\ldots,n_k} \subset B_{n_1,n_2,\ldots,n_k}$ such that $h(P_{n_1,n_2,\ldots,n_k})=3 \cdot \prod_{i=1}^{k-1} n_i-2$ and \linebreak$l(P)=n_k-1$.

$P_{n_1,n_2,\ldots,n_k}$ [First Layer] $= (0,\hspace{2.5mm}0,\hspace{2.5mm}0,\hspace{2.5mm}\ldots,\hspace{2.5mm}0,\hspace{2.5mm}0)$-$(0,\hspace{2.5mm}0,\hspace{2.5mm}0,\hspace{2.5mm}\ldots,\hspace{2.5mm}0,\hspace{2.5mm}n_k-1)$-\linebreak$\bigg(\frac{1}{2},\hspace{1.5mm}0,\hspace{1.5mm}0,\hspace{1.5mm}\ldots,\hspace{1.5mm}0,\hspace{1.5mm}n_k-1-\sqrt{(n_k-1)^2-\frac{1}{4}})\bigg)$-$(1,\hspace{1.5mm}0,\hspace{1.5mm}0,\hspace{1.5mm}\ldots,\hspace{1.5mm}0,\hspace{1.5mm}n_k-1)$-$(1,\hspace{1.5mm}0,\hspace{1.5mm}0,\hspace{1.5mm}\ldots,\hspace{1.5mm}0,\hspace{1.5mm}0)$-\linebreak$\bigg(\frac{3}{2},\hspace{2.7mm}0,\hspace{2.7mm}0,\hspace{2.7mm}\ldots,\hspace{2.7mm}0,\hspace{2.7mm}\sqrt{((n_k-1)^2-\frac{1}{4}})\bigg)$-$(2,\hspace{2.7mm}0,\hspace{2.7mm}0,\hspace{2.7mm}\ldots,\hspace{2.7mm}0,\hspace{2.7mm}0)$-$(2,\hspace{2.7mm}0,\hspace{2.7mm}0,\hspace{2.7mm}\ldots,\hspace{2.7mm}0,\hspace{2.7mm}n_k-1)$-\linebreak$\bigg(\frac{5}{2},\hspace{1.5mm}0,\hspace{1.5mm}0,\hspace{1.5mm}\ldots,\hspace{1.5mm}0,\hspace{1.5mm}n_k-1-\sqrt{(n_k-1)^2-\frac{1}{4}}\bigg)$-$(3,\hspace{1.5mm}0,\hspace{1.5mm}0,\hspace{1.5mm}\ldots,\hspace{1.5mm}0,\hspace{1.5mm}n_k-1)$-$(3,\hspace{1.5mm}0,\hspace{1.5mm}0,\hspace{1.5mm}\ldots,\hspace{1.5mm}0,\hspace{1.5mm}0)$-\linebreak$\bigg(\frac{7}{2},0,0,\ldots,0,\sqrt{(n_k-1)^2-\frac{1}{4}}\bigg)$-$\hspace{0.5mm}\ldots$ and so forth, until either $(n_1-1,\hspace{1mm}0,\hspace{1mm}0,\hspace{1mm}\ldots,\hspace{1mm}0,\hspace{1mm}n_k-1)$ or $(n_1-1,\hspace{1mm}0,\hspace{1mm}0,\hspace{1mm}\ldots,\hspace{1mm}0,\hspace{1mm}0)$ is reached last (depending on whether $n_1$ is odd or even, respectively).

Now, let us assume that $n_1$ is even so that the last link to cover the first layer is \linebreak$(n_1-1,\hspace{1mm}0,\hspace{1mm}0,\hspace{1mm}\ldots,\hspace{1mm}0,\hspace{1mm}n_k-1)$-$(n_1-1,\hspace{1mm}0,\hspace{1mm}0,\hspace{1mm}\ldots,\hspace{1mm}0,\hspace{1mm}0)$. Otherwise, $n_1$ has to be odd and we fall in the opposite case, where two flipped bridges will alternately be switched moving from layer to layer.

Then, by assuming for simplicity $n_1 : n_1=2 \cdot m$, where $m \in \mathbb{N}-\{0,1\}$, we have \linebreak
$P_{n_1,n_2,\ldots,n_k}$ [First Bridge]$\hspace{2mm} =\bigg(n_1-1,\hspace{1mm}\frac{1}{2},\hspace{1mm}0,\hspace{1mm}\ldots,\hspace{1mm}0,\hspace{1mm}\sqrt{(n_k-1)^2-\frac{1}{4}}\bigg)$-$(n_1-1,\hspace{1mm}1,\hspace{1mm}0,\hspace{1mm}\ldots,\hspace{1mm}0,\hspace{1mm}0)$.\linebreak
$P_{n_1,\hspace{1mm}n_2,\hspace{1mm}\ldots,\hspace{1mm}n_k}$ [Second Layer]\hspace{2mm}$=(n_1-1,\hspace{2mm}1,\hspace{2mm}0,\hspace{2mm}\ldots,\hspace{2mm}0,\hspace{2mm}n_k-1)$-$\bigg(n_1-1-\frac{1}{2},\hspace{2mm}1,\hspace{2mm}0,\hspace{2mm}\ldots,\hspace{2mm} \linebreak 0,\hspace{2mm}n_k-1-\sqrt{(n_k-1)^2-\frac{1}{4}}\bigg)$-$(n_1-2,\hspace{2mm}1,\hspace{2mm}0,\hspace{2mm}\ldots,\hspace{2mm}0,\hspace{2mm}n_k-1)$-$(n_1-2,\hspace{2mm}1,\hspace{2mm}0,\hspace{2mm}\ldots,\hspace{2mm}0,\hspace{2mm}0)$-\linebreak$\bigg(n_1-2-\frac{1}{2},\hspace{1mm}1,\hspace{1mm}0,\hspace{1mm}\ldots,\hspace{1mm}0,\hspace{1mm}\sqrt{(n_k-1)^2-\frac{1}{4}}\bigg)$-$(n_1-3,\hspace{1mm}1,\hspace{1mm}0,\hspace{1mm}\ldots,\hspace{1mm}0,\hspace{1mm}0)$-$(n_1-3,\hspace{1mm}1,\hspace{1mm}0,\hspace{1mm}\ldots,\hspace{1mm}0,\hspace{1mm}n_k-1)-\ldots$ and so on. In this way, we finally reach $(0,0,1,\ldots,0,0)$, the last point of the second layer.

$P_{n_1,\hspace{1mm}n_2,\hspace{1mm}\ldots,\hspace{1mm}n_k}$ [Second Bridge] $=(0,\hspace{1mm}0,\hspace{1mm}\frac{3}{2},\hspace{1mm}\ldots,\hspace{1mm}0,\hspace{1mm}\sqrt{(n_k-1)^2-\frac{1}{4}})$-$(0,\hspace{1mm}0,\hspace{1mm}2,\hspace{1mm}\ldots,\hspace{1mm}0,\hspace{1mm}0)$.

We repeat the same pattern until we reach the endpoint, which is the last visited point among all the elements of the set $V$ defined as\vspace{-2mm}
$$V \coloneqq \big\{(x_1,x_2,\ldots,x_k) : x_1 \in \{0,n_1-1\}  \wedge x_2 \in \{0, n_2-1\} \wedge \cdots \wedge x_k \in \{0, n_k-1\}\big\}.$$

Therefore, $\forall(n_1,n_2,\ldots,n_k) : n_1 \leq n_2 \leq \cdots \leq n_k \in \mathbb{N}-\{0,1\}$, $\exists P_{n_1,n_2,\ldots,n_k}$, an inside the MAABB uncrossing covering path, such that $h(P_{n_1,n_2,\ldots,n_k} )<\prod_{i=1}^{k} n_i \wedge l(P_{n_1,n_2,\ldots,n_k})= 
n_k-1$. Moreover, the M$\bf{\Lambda}$I-algorithm let $P_{n_1,n_2,\ldots,n_k}$ join an increasing number of points with every triplet of consecutive edges as $n_k$ grows, and Equations (\ref{eq4})\&(\ref{eq5}) show how its efficiency becomes absolute when $n_k$ approaches infinity. Then from Equation (\ref{eq2}), it follows that
\begin{equation} \label{eq4}
\frac{h(P_{n_1,n_2,\ldots,n_k})}{\prod_{i=1}^{k} n_i}=\frac{3}{n_k}-\frac{2}{\prod_{i=1}^{k} n_i}.
\end{equation}

Hence,
\begin{equation} \label{eq5}
\lim_{n_k\to +\infty} \frac{h(P_{n_1,n_2,\ldots,n_k})}{\prod_{i=1}^{k} n_i}=0.
\end{equation}

\begin{figure}[H]
\begin{center}
\includegraphics[scale=0.85]{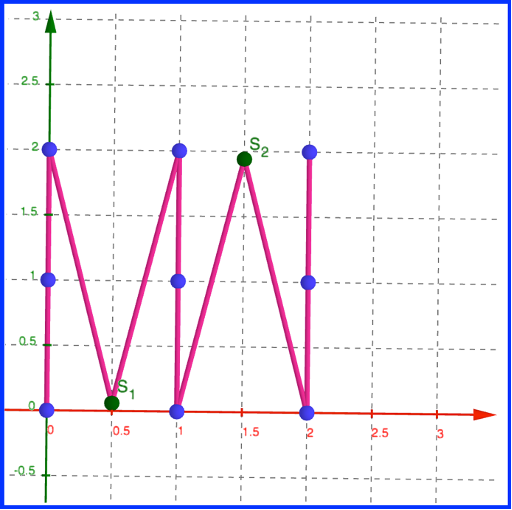}
\end{center}
\caption{First layer of the covering path $P_{3,3,3}$, edges $1$ to $7$: all of them belong to the length class $2$. The two Steiner points (in green) are $S_1 \equiv \left(\frac{1}{2},2-\frac{\sqrt{15}}{2}\right) \cong (0.5, 0.063508)$ and $S_2 \equiv \left(\frac{3}{2},\frac{\sqrt{15}}{2}\right) \cong (1.5, 1.936492)$ \protect\cite{8}.}
\label{fig:Figure_1c}
\end{figure}

\begin{figure}[H]
\begin{center}
\includegraphics{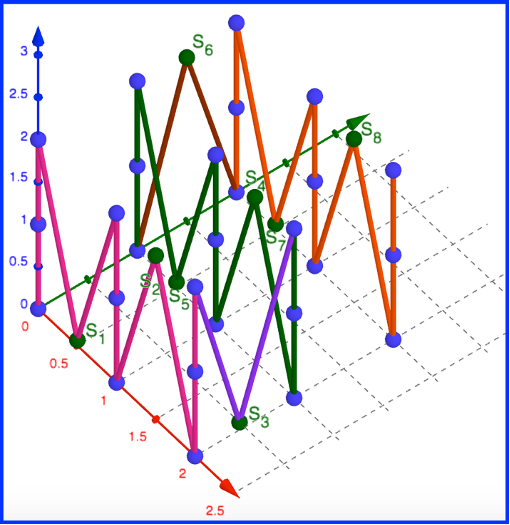}
\end{center}
\caption{The whole directed covering path $P_{3,3,3}$, edges $1$ to $25$: all of them belong to the length class $2$ \protect\cite{8}.}
\label{fig:Figure_2c}
\end{figure}

As a couple of examples, we can take $P_{3,3} \subset \hat{B}_{3,3}$ and $P_{3,3,3} \subset \hat{B}_{3,3,3} \hspace{1mm}$, where \linebreak $P_{3,3}=(0,0)$-$(0,2)$-$\left(\frac{1}{2},2-\frac{\sqrt{15}}{2}\right)$-$(1,2)$-$(1,0)$-$\left(\frac{3}{2},\frac{\sqrt{15}}{2}\right)$-$(2,0)$-$(2,2)$ and $P_{3,3,3}=(0,0,0)$-$(0,0,2)$-\linebreak $\left(\frac{1}{2},\hspace{0.4mm}0,\hspace{0.4mm}2-\frac{\sqrt{15}}{2}\right)$-$(1,\hspace{0.4mm}0,\hspace{0.4mm}2)$-$(1,\hspace{0.4mm}0,\hspace{0.4mm}0)$-$\left(\frac{3}{2},\hspace{0.4mm}0,\hspace{0.4mm}\frac{\sqrt{15}}{2}\right)$-$(2,\hspace{0.4mm}0,\hspace{0.4mm}0)$-$(2,\hspace{0.4mm}0,\hspace{0.4mm}2)$-$\left(2,\hspace{0.4mm}\frac{1}{2},\hspace{0.4mm}2-\frac{\sqrt{15}}{2}\right)$-$(2,\hspace{0.4mm}1,\hspace{0.4mm}2)$-$(2,\hspace{0.4mm}1,\hspace{0.4mm}0)$-\linebreak$\left(\frac{3}{2},\hspace{0.75mm}1,\hspace{0.75mm}\frac{\sqrt{15}}{2}\right)$-$(1,\hspace{0.75mm}1,\hspace{0.75mm}0)$-$(1,\hspace{0.75mm}1,\hspace{0.75mm}2)$-$\left(\frac{1}{2},\hspace{0.75mm}1,\hspace{0.75mm}2-\frac{\sqrt{15}}{2}\right)$-$(0,\hspace{0.75mm}1,\hspace{0.75mm}2)$-$(0,\hspace{0.75mm}1,\hspace{0.75mm}0)$-$\left(0,\hspace{0.75mm}\frac{3}{2},\hspace{0.75mm}\frac{\sqrt{15}}{2}\right)$-$(0,\hspace{0.75mm}2,\hspace{0.75mm}0)$-$(0,\hspace{0.75mm}2,\hspace{0.75mm}2)$-\linebreak$\left(\frac{1}{2},\hspace{0.5mm}2,\hspace{0.5mm}2-\frac{\sqrt{15}}{2}\right)$-$(1,\hspace{0.5mm}2,\hspace{0.5mm}2)$-$(1,\hspace{0.5mm}2,\hspace{0.5mm}0)$-$\left(\frac{3}{2},\hspace{0.5mm}2,\hspace{0.5mm}\frac{\sqrt{15}}{2}\right)$-$(2,\hspace{0.5mm}2,\hspace{0.5mm}0)$-$(2,\hspace{0.5mm}2,\hspace{0.5mm}2)$, as shown in Figures \ref{fig:Figure_1c}\&\ref{fig:Figure_2c}.

\begin{figure}[H]
\begin{center}
\includegraphics{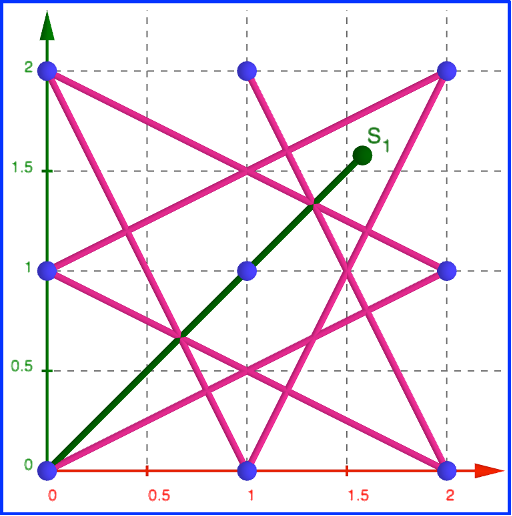}
\end{center}
\caption{The self-intersecting covering path $M_{3,3}$ consists of $8$ edges, and all of them belong to the length class $\sqrt{5}$ \protect\cite{8}.}
\label{fig:Figure_3c}
\end{figure}

Referring to $G_{3,3}$ and $G_{3,3,3}$ (as above), it is possible to show also the existence of self-intersecting covering paths, $M_{3,3} \subset \hat{B}_{3,3}$ and $M_{3,3,3} \subset \hat{B}_{3,3,3}$, such that their link length is $\prod_{i=1}^{k} 3^i -1$ and $l(M_{3,3})=l(M_{3,3,3})=\sqrt{5}>n_ k-1$. In particular, as shown in Figure \ref{fig:Figure_3c}, $M_{3,3}\hspace{1mm}=\hspace{1mm}(1,\hspace{1.5mm}2)$-$(2,\hspace{1.5mm}0)$-$(0,\hspace{1.5mm}1)$-$(2,\hspace{1.5mm}2)$-$(1,\hspace{1.5mm}0)$-$(0,\hspace{1.5mm}2)$-$(2,\hspace{1.5mm}1)$-$(0,\hspace{1.5mm}0)$-$\left(\sqrt{\frac{5}{2}},\hspace{1.5mm}\sqrt{\frac{5}{2}}\right)$, while \linebreak$M_{3,3,3}=(2,\hspace{0.4mm}0,\hspace{0.3mm}0)$-$(0,\hspace{0.3mm}1,\hspace{0.3mm}0)$-$(2,\hspace{0.3mm}2,\hspace{0.3mm}0)$-$(1,\hspace{0.3mm}0,\hspace{0.3mm}0)$-$(0,\hspace{0.3mm}2,\hspace{0.3mm}0)$-$(1,\hspace{0.3mm}2,\hspace{0.3mm}2)$-$(0,\hspace{0.3mm}0,\hspace{0.3mm}2)$-$(2,\hspace{0.3mm}1,\hspace{0.3mm}2)$-$(0,\hspace{0.3mm}2,\hspace{0.3mm}2)$-$(1,\hspace{0.3mm}0,\hspace{0.3mm}2)$-\linebreak$(2,\hspace{0.33mm}2,\hspace{0.33mm}2)$-$(0,\hspace{0.33mm}1,\hspace{0.33mm}2)$-$(2,\hspace{0.33mm}0,\hspace{0.33mm}2)$-$(2,\hspace{0.33mm}2,\hspace{0.33mm}1)$-$(0,\hspace{0.33mm}1,\hspace{0.33mm}1)$-$(2,\hspace{0.33mm}0,\hspace{0.33mm}1)$-$(1,\hspace{0.33mm}2,\hspace{0.33mm}1)$-$(0,\hspace{0.33mm}0,\hspace{0.33mm}1)$-$(2,\hspace{0.33mm}1,\hspace{0.33mm}1)$-$(0,\hspace{0.33mm}2,\hspace{0.33mm}1)$-$(1,\hspace{0.33mm}0,\hspace{0.33mm}1)$-\linebreak$(1,\hspace{-0.1mm}2,\hspace{-0.1mm}0)$-$(1,\hspace{-0.1mm}1,\hspace{-0.1mm}2)$-$(2,\hspace{-0.1mm}1,\hspace{-0.1mm}0)$-$(0,\hspace{-0.1mm}0,\hspace{-0.1mm}0)$-$\left(\sqrt{\frac{5}{2}},\hspace{-0.1mm}\sqrt{\frac{5}{2}},\hspace{-0.1mm}0\right)$-$\Big(\frac{\sqrt{10}}{12} \hspace{-0.1mm}\cdot \hspace{-0.1mm} \big(6-\sqrt{24-3 \cdot \sqrt{10}}\big),\hspace{-0.1mm}\frac{\sqrt{10}}{12} \hspace{-0.1mm} \cdot \hspace{-0.1mm} \big(6-\sqrt{24-3 \cdot \sqrt{10}} \big),\linebreak \frac{1}{2} \cdot \sqrt{\frac{5}{3} \cdot (4+\sqrt{10}})\Big)$ (see \cite{3}, proof of Theorem 2, for details).

\begin{figure}[H]
\begin{center}
\includegraphics{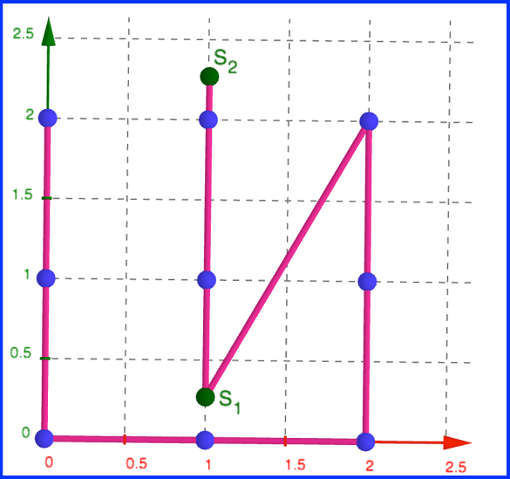}
\end{center}
\caption{The covering path $\check{P}_{3,3}=(0,2)$-$(0,0)$-$(2,0)$-$(2,2)$-$(1,2-\sqrt{3})$-$(1,4-\sqrt{3})$ \protect\cite{8}.}
\label{fig:Figure_4c}
\end{figure}

\begin{figure}[H]
\begin{center}
\includegraphics{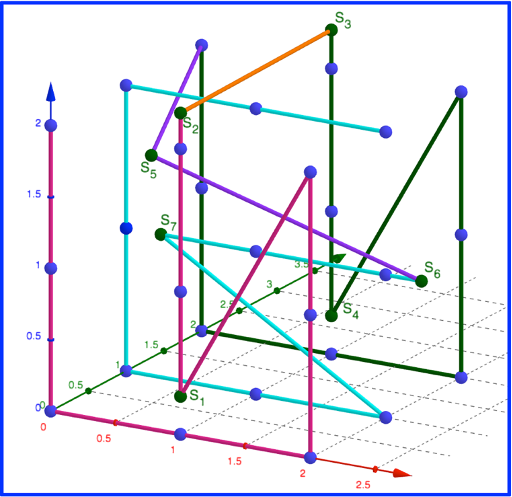}
\end{center}
\caption{The (uncrossing) covering path $\check{P}_{3,3,3}=(0,2,0)$-$(0,0,0)$-$(2,0,0)$-$(2,2,0)$-$(1,2-\sqrt{3},0)$-$(1,4-\sqrt{3},0)$-$(1,4-\sqrt{3},2)$-$(1,2-\sqrt{3},2)$- \linebreak$(2,2,2)$-$(2,0,2)$-$(0,0,2)$-$(0,2,2)$-$\left(\frac{7}{10},\frac{33 \cdot \sqrt{3}-37+\sqrt{2 \cdot \left(2541 \cdot \sqrt{3}-4247 \right)}}{20},\frac{33 \cdot \sqrt{3}-37-\sqrt{2 \cdot \left(2541 \cdot \sqrt{3}-4247 \right)}}{20}\right)$- \linebreak$(4-\sqrt{3},1,1)$-$(2-\sqrt{3},1,1)$-$(2,0,1)$-$(0,0,1)$-$(0,2,1)$-$(2,2,1)$ extends $\check{P}_{3,3}$ to $3$ dimensions \protect\cite{8}.}
\label{fig:Figure_5c}
\end{figure}

On the other hand, considering all the grids $G_{3,3,\ldots,3}(k) : k \geq 2$, if we loosen the constraint on the minimum axis-aligned bounding box to the regular one (or at least to $\check{B}_{3,3,\ldots,3} \coloneqq [0,2] \times \left[0,4-\sqrt{3}\right] \times \left[0,4-\sqrt{3}\right] \times [0,2]^{k-3}$), we can easily show the existence of uncrossing covering paths with at most $h(\check{P}_{3,3,\ldots,3})=\ceil*{20 \cdot 3^{k-3}}-2$ links, all belonging to the length class $l(\check{P}_{3,3,\ldots,3})=2$ (see Figures \ref{fig:Figure_4c}\&\ref{fig:Figure_5c}).

Consider the uncrossing covering path $\check{P}_{3,3,3}$ (see Figure \ref{fig:Figure_5c} where $x_{S_5}=0.7$ and the bigger of the two values of $y(x_{S_5})$ have been selected). The Steiner point $S_5 \equiv \left(x_{S_5},y_{S_5},z_{S_5}\right)$ can be arbitrarily chosen among many of the solutions provided by the intersection of the two spheres of radius $2$ units which are centered in $A \equiv (0,2,2)$ and $B \equiv (4-\sqrt{3},1,1)$. Thus,
\begin{equation} \label{eq6}
 \begin{cases}
\hspace{0.5 mm} x^2+(y-2)^2+(z-2)^2=2^2 \\
\hspace{0.5 mm}  \left(x-4+\sqrt{3}\right)^2+(y-1)^2+(z-1)^2=2^2 \\
 \end{cases}.
\end{equation}

For this purpose, it is necessary to preliminary point out that the two conditions $y_{S_5} \in \left[\frac{3}{2}-\sqrt{\frac{466 \cdot \sqrt{3}-333}{249}},4-\sqrt{3} \right]-\{2\}$ and $z_{S_5}\leq 2$ should be satisfied, and they are a sufficient restriction to ensure that any constrained solution of (\ref{eq7}) returns a valid $S_5$, avoiding any self-intersecting risk for $\check{P}_{3,3,3}$ (see Relations (\ref{eq7}) to (\ref{eq9}) below), so we have
\begin{flalign} \label{eq7}
\begin{aligned} \raisetag{1ex}
y_{S_5}\hspace{-1mm}=\hspace{-1mm}\frac{\sqrt{(32 \cdot \sqrt{3}-84) \cdot x^2+(432-212 \cdot \sqrt{3}) \cdot x+336 \cdot \sqrt{3}-601}+ (8-2 \cdot\sqrt{3}) \cdot x+8 \cdot \sqrt{3}-13}{4} \wedge \\[4pt]
z_{S_5}\hspace{-1mm}=\hspace{-1mm}\frac{-\sqrt{(32 \cdot \sqrt{3}-84) \cdot x^2+(432-212 \cdot \sqrt{3}) \cdot x+336 \cdot \sqrt{3}-601}+ (8-2 \cdot\sqrt{3}) \cdot x+8 \cdot \sqrt{3}-13}{4} \hspace{-0.1mm}.
\end{aligned}
\end{flalign}

More specifically, the circumference by (\ref{eq6}) under the pair of previously stated conditions assures the existence of (at least) one possible Steiner point $S_5$, which is well-defined assuming $S_5 \equiv \overline{S}_5 \vee \left(y_{S_5} \geq z_{S_5}  \wedge S_5 \not\equiv \overline{\overline{S}}_5 \right)$ (as above), for every
\begin{equation} \label{eq8}
x \in \left[2-\frac{\sqrt{3}}{2}-\sqrt{\frac{87+128 \cdot \sqrt{3}}{498}}, \hspace{1mm} \frac{65 \cdot \left(2-\sqrt{3}\right)-\sqrt{916 \cdot \sqrt{3}-1549}}{8 \cdot \left(5-2 \cdot \sqrt{3} \right)} \right] \approx \left[0.346647, 0.918696 \right].
\end{equation}

\begin{figure}[H]
\begin{center}
\includegraphics[scale=1.15]{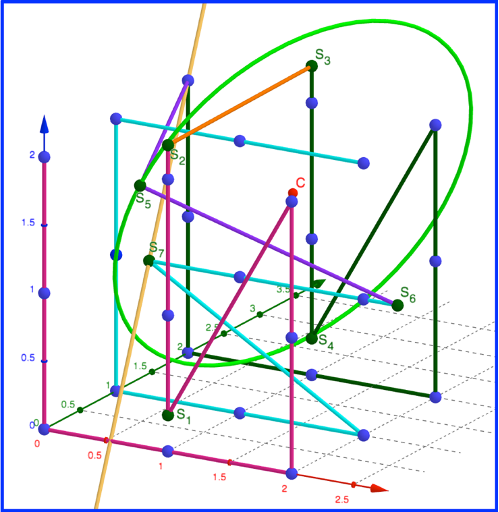}
\end{center}
\caption{Graphical proof that no collision will occur in the process of moving $S_5$ along the allowed arc $\left(2>y_{S_5} \geq z_{S_5} \right)$ of the circle given by Equation (\ref{eq6}), as shown by comparison of (\ref{eq6}) with (\ref{eq9}) (i.e., observing that $2>r \approx 1.9715304811$). The radius of the circle centered in $C \equiv \left(\frac{4-\sqrt{3}}{2},\frac{3}{2},\frac{3}{2} \right)$ is $a=\frac{\sqrt{8 \cdot \sqrt{3}-5}}{2} \approx 1.4879857577$, while a collision between the $13$th and the $15$th link would only occur for $y_{S_5}=\frac{161+2 \cdot \sqrt{3}-\sqrt{20092 \cdot \sqrt{3}-17383}}{104} < 1 = z_{S_5}$ \protect\cite{8}.}
\label{fig:Figure_6c}
\end{figure}

Forasmuch as $\overline{\overline{S}}_5 \equiv \left(\frac{53 \cdot \sqrt{3}-108+ \sqrt{208 \cdot \sqrt{3}-313}}{8 \cdot \left(2 \cdot \sqrt{3}-5 \right)},\hspace{2mm} 2, \hspace{2mm}\frac{151-2 \cdot \sqrt{3}-\sqrt{20092 \cdot \sqrt{3}-17383}}{104} \right)$ cannot be accepted as a solution in order to generate a covering path which is not self-intersecting (i.e., $y_{S_5} \neq 2$ by one of the stated conditions), let us replace the invalid Steiner point $\overline{\overline{S}}_5$ with $\overline{S}_5 \equiv \left(\frac{53 \cdot \sqrt{3}-108+ \sqrt{208 \cdot \sqrt{3}-313}}{8 \cdot \left(2 \cdot \sqrt{3}-5 \right) }, \hspace{2mm}\frac{151-2 \cdot \sqrt{3}-\sqrt{20092 \cdot \sqrt{3}-17383}}{104}, \hspace{2mm} 2 \right)$. Finally, we observe that there are not any other critical values, because (see Figure \ref{fig:Figure_6c}, looking for any possible collision between the twelfth link of $\check{P}_{3,3,3}$ and the light blue layer in the middle)
\begin{equation} \label{eq9}
\begin{aligned}
\hspace{-5cm} \begin{cases}
x^2+(y-2)^2+(z-2)^2=r^2 \\
\left(x-4+\sqrt{3} \right)^2+(y-1)^2+(z-1)^2=r^2
\end{cases}  \cap  \quad \hspace{1.8mm}
\begin{cases}
\frac{x}{2-\sqrt{3}} =2-y \\
y=z
\end{cases}  \Rightarrow \hspace{27.5mm}
\\
\begin{cases}
x^2+(y-2)^2+(z-2)^2=r^2 \\
\left(x-4+\sqrt{3}\right)^2+(y-1)^2+(z-1)^2=r^2 \\
\frac{x}{2-\sqrt{3}}=2-y \\
y=z
\end{cases}  \Rightarrow \quad
\begin{cases}
r=\frac{\sqrt{94725-21288 \cdot \sqrt{3}}}{122} \\
x=\frac{192-85 \cdot \sqrt{3}}{122} \\
y=\frac{115-22 \cdot \sqrt{3}}{122}  \\
z=\frac{115-22 \cdot \sqrt{3}}{122}
\end{cases} \Rightarrow \quad r < 2.
\end{aligned}
\end{equation}

Let us extend the results to higher dimensions, we then assume $k \geq 3$.

From Equations (\ref{eq7})\&(\ref{eq8}), since it is possible to set an unlimited number of layer bridges by choosing distinct Steiner points as $S_5$, we can straightforwardly derive uncrossing covering paths $\check{P}_{3,3,\ldots,3}$ with $18 \cdot 3^{k-3}+\sum_{j=0}^{k-4} 3^j \cdot 4= 20 \cdot 3^{k-3}-2$ edges if we simply lift, by iteration, the given $3$D pattern to the next dimension (in a very similar way as we previously shown in Figure \ref{fig:Figure_2c} for the M$\Lambda$I-algorithm). Furthermore, referring to the same AABB (i.e., $\check{B}_{3,3,\ldots,3}$) and pattern by Figures \ref{fig:Figure_5c}\&\ref{fig:Figure_6c}, it is possible to save one more line for any dimensional bridge if we switch from uncrossing covering paths to covering trails, spending a total of $18 \cdot 3^{k-3}+\sum_{j=1}^{k-3} 3^j=  \frac{13 \cdot 3^{k-2} - 3}{2}$ links (at most) to join all the given $3^k$ points, and this consideration would suffice to reveal how many unexplored optimization problems could be studied starting from a different set of preliminary assumptions, rather than the four fundamental constraints stated in Section \ref{sec:Intr}.

\begin{figure}[H]
\begin{center}
\includegraphics{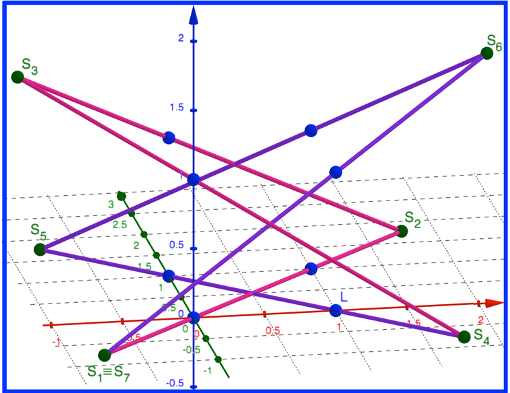}
\end{center}
\caption{A smart covering cycle for $G_{2,2,2}$. $F_{2,2,2}$ has a link length of $6$, is constrained by $\left[-\frac{1+\sqrt{13}}{4},\frac{5+\sqrt{13}}{4}  \right] \times \left[-\frac{1+\sqrt{13}}{6},\frac{7+\sqrt{13}}{6} \right] \times \left[0,\frac{3+\sqrt{13}}{4} \right] \not\subseteq $ RAABB, and is characterized by
$l(F_{2,2,2})=\frac{\sqrt{2} \cdot \left(4+\sqrt{13}\right)}{3}$ (see Definitions \ref{def3} to \ref{def6}) 
 \protect\cite{8}.}
\label{fig:Figure_7c}
\end{figure}

\begin{figure}[H]
\begin{center}
\includegraphics{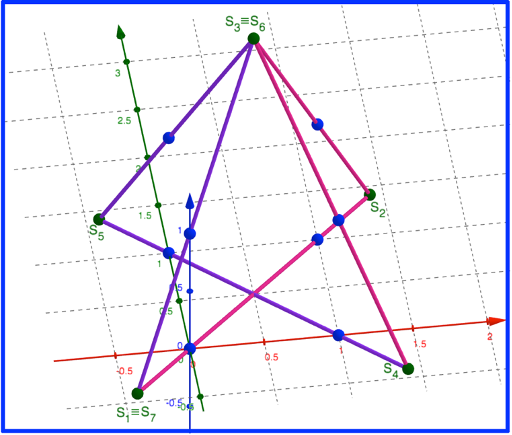}
\end{center}
\caption{A covering circuit with optimal link length $h(F'_{2,2,2})=6$. $F'_{2,2,2}$ is inside the AABB $\left[1-\sqrt{2},\sqrt{2} \right] \times \left[1-\sqrt{2},\sqrt{2} \right] \times \left[0,2 \cdot \sqrt{3}-\sqrt{\frac{3}{2}} \right]\not\subseteq$ RAABB and is characterized by $l(F'_{2,2,2})=4-\sqrt{2}$ (see Definitions \ref{def3} to \ref{def6}) \protect\cite{8}.}
\label{fig:Figure_8c}
\end{figure}

We could also go further and decide (as a random example) to cover the grid with an unconstrained circuit of minimum link length, disregarding the RAABB,
so we may accept $F_{2,\hspace{0.2mm}2,\hspace{0.2mm}2}=\left(-\frac{1+\sqrt{13}}{6},\hspace{0.2mm}-\frac{1+\sqrt{13}}{6},\hspace{0.2mm}0\right)$-$\left(\frac{7+\sqrt{13}}{6},\hspace{0.2mm}\frac{7+\sqrt{13}}{6},\hspace{0.2mm}0 \right)$-$\left(-\frac{1+\sqrt{13}}{4},\hspace{0.2mm}\frac{1}{2},\hspace{0.2mm}\frac{3+\sqrt{13}}{4}\right)$-$\left(\frac{7+\sqrt{13}}{6},\hspace{0.2mm}-\frac{1+\sqrt{13}}{6},\hspace{0.2mm}0 \right)$- \linebreak$\left(-\frac{1+\sqrt{13}}{6},\hspace{0.7mm}\frac{7+\sqrt{13}}{6},\hspace{0.7mm}0\right)$-$\left(\frac{5+\sqrt{13}}{4},\hspace{0.7mm}\frac{1}{2},\hspace{0.7mm}\frac{3+\sqrt{13}}{4} \right)$-$\left(-\frac{1+\sqrt{13}}{6},\hspace{0.7mm}-\frac{1+\sqrt{13}}{6},\hspace{0.7mm}0\right)$ and $F'_{2,\hspace{0.7mm}2,\hspace{0.7mm}2}=\left(1-\sqrt{2},\hspace{0.7mm}1-\sqrt{2},\hspace{0.7mm}0\right)$- \linebreak$\left(\sqrt{2},\hspace{2.1mm} \sqrt{2},\hspace{2.1mm} 0 \right)$-$\left(\frac{1}{2},\hspace{2.1mm}\frac{1}{2},\hspace{2.1mm}2 \cdot \sqrt{3}-\sqrt{\frac{3}{2}}\right)$-$\left(\sqrt{2},\hspace{2.1mm}1-\sqrt{2},\hspace{2.1mm} 0 \right)$-$\left(1-\sqrt{2},\hspace{2.2mm}\sqrt{2},\hspace{2.1mm}0 \right)$-$\left(\frac{1}{2},\hspace{2.2mm}\frac{1}{2},\hspace{2.2mm}2 \cdot \sqrt{3}-\sqrt{\frac{3}{2}} \right)$- \linebreak$\left(1-\sqrt{2}, \hspace{2mm} 1-\sqrt{2}, \hspace{2mm} 0 \right)$ as a pair of valid solutions (see Figures \ref{fig:Figure_7c}\&\ref{fig:Figure_8c}), while we may reject $(0,0,0)$-$(1,0,0)$-$(1,1,0)$-$(0,1,0)$-$(0,1,1)$-$(1,1,1)$-$(1,0,1)$-$(0,0,1)$ for two different reasons.

Lastly, it is even possible to approach some optimization problems of path covering focusing on the AABB volume, such as finding the uncrossing covering path of minimal link length which subtends the minimum volume AABB.

Let $G_{2,2,2}$ be given and assume $x_{S_1 }>0$. We could begin with the set of self-intersecting covering paths $\overline{P}_{2,\hspace{0.4mm}2,\hspace{0.4mm}2}\left(x_{S_1} \right)=(0,\hspace{0.4mm}1,\hspace{0.4mm}0)$-$(0,\hspace{0.4mm}0,\hspace{0.4mm}0)$-$\left(x_{S_1},\hspace{0.4mm}0,\hspace{0.4mm}x_{S_1} \right)$-$\left(\frac{1}{2},\hspace{0.4mm}y_{S_2}\left(x_{S_1} \right),\hspace{0.4mm}\frac{1}{2}\right)$-$\left(1-x_{S_1},\hspace{0.4mm}0,\hspace{0.4mm}x_{S_1} \right)$-\linebreak$(1,\hspace{0.4mm}0,\hspace{0.4mm}0)$-$(1,\hspace{0.4mm}1,\hspace{0.4mm}0)$, and then determine the value of $x_{S_1}$ that minimizes the volume of the AABB from
\begin{equation} \label{eq10}
\begin{cases}
\left(1-x_{S_1}\right) \cdot y_{S_2}+x_{S_1}=\frac{1}{2} \\[6pt]
\frac{\partial{\left({{\left(2 \cdot x_{S_1}-1 \right)} \cdot y_{S_2} \cdot x_{S_1}}\right)}}{\partial{x_{S_1}}}=0 \\
x_{S_1}>1
\end{cases}.
\end{equation}

Hence,
\begin{equation} \label{eq11}
\begin{aligned}
\begin{cases} y_{S_2}=\frac{1 - 2 \cdot x_{S_1}}{2 \cdot \left(1 - x_{S_1} \right)} \\
{x_{S_1}}^3-2 \cdot {x_{S_1}}^2+x_{S_1}=\frac{1}{8} \\
x_{S_1}>1
\end{cases}
\Rightarrow \quad
\begin{cases}
x_{S_1}=\frac{3 + \sqrt{5}}{4} \\
y_{S_2}=\frac{3 + \sqrt{5}}{2}
\end{cases}.
\end{aligned}
\end{equation}

Now, let $\varphi \coloneqq \frac{1 + \sqrt{5}}{2}$ (i.e., the well-known \textit{golden ratio}).

It follows that $x_{S_1}=\frac{1 + \varphi}{2} \hspace{1mm} \wedge \hspace{1mm} y_{S_2}=1+\varphi$.

\begin{figure}[H]
\begin{center}
\includegraphics{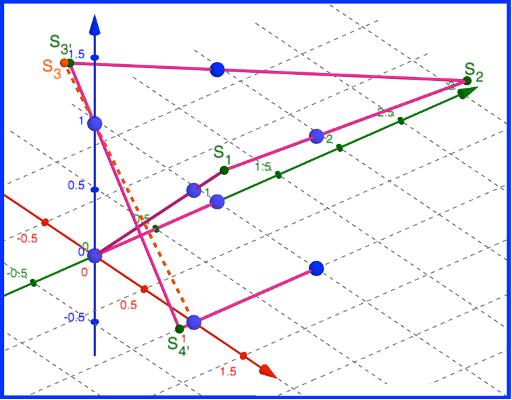}
\end{center}
\caption{The self-crossing covering path $\overline{P}_{2,2,2} \left( \frac{1 + \varphi}{2}\right)$ and the uncrossing covering path $\overline{\overline{P}}_{2,2,2}(\varepsilon)$, where the (Euclidean) distance between the point $(1,0,0)$ and the Steiner point $S_{4'} \in \overline{\overline{P}}_{2,2,2} (\varepsilon)$ is given by $-y_{S_{4'}}=-\frac{4 \cdot \varphi \cdot \varepsilon}{1 - \varphi + 2 \cdot \varepsilon}$ \protect\cite{8}.}
\label{fig:Figure_9c}
\end{figure}

Thus, for any $0<\varepsilon<\frac{\varphi - 1}{2}$, we can easily derive the uncrossing covering path $\overline{\overline{P}}_{2,2,2} (\varepsilon)=(0,\hspace{0.35mm}1,\hspace{0.35mm}0)$-$(0,\hspace{0.35mm}0,\hspace{0.35mm}0)$-$\left(\frac{1 + \varphi}{2},\hspace{0.35mm}0,\hspace{0.35mm}\frac{1 + \varphi}{2} \right)$-$\left(\frac{1}{2},\hspace{0.35mm}1+\varphi,\hspace{0.35mm}\frac{1}{2} \right)$-$\left(\frac{1 - \varphi}{2}+\varepsilon,\hspace{0.35mm}2 \cdot \varphi \cdot \varepsilon,\hspace{0.35mm}\frac{1 + \varphi}{2}-\varepsilon \right)$-$\left(1,\hspace{0.35mm}\frac{4 \cdot \varphi \cdot \varepsilon}{1 - \varphi + 2 \cdot \varepsilon},\hspace{0.35mm}0 \right)$-$(1,\hspace{0.35mm}1,\hspace{0.35mm}0)$ from $\overline{P}_{2,2,2} \left(\frac{1 + \varphi}{2} \right)=(0,1,0)$-$(0,0,0)$-$\left(\frac{1 + \varphi}{2},0,\frac{1 + \varphi}{2} \right)$-$\left(\frac{1}{2},1+\varphi,\frac{1}{2} \right)$-$\left(\frac{1 - \varphi}{2},0,\frac{1 + \varphi}{2} \right)$-$(1,0,0)$-$(1,1,0)$. Consequently, $\lim_{\varepsilon \to 0}\overline{\overline{P}}_{2,2,2} (\varepsilon)=\overline{P}_{2,2,2} \left(\frac{1+\varphi}{2} \right)$ (see Figure \ref{fig:Figure_9c}).

Therefore, $\overline{\overline{P}}_{2,2,2}(\varepsilon)$ shows the existence of minimal uncrossing covering paths which are entirely contained in $\left[\frac{1 - \varphi}{2}+ \varepsilon, \frac{1 + \varphi}{2} \right] \times \left[\frac{4 \cdot \varphi \cdot \varepsilon}{1 - \varphi+ 2 \cdot \varepsilon},1+\varphi \right] \times \left[0,\frac{1 + \varphi}{2} \right] \subset [-0.309016, 1.309017]  \times [-10.47214 \cdot \varepsilon, 2.618034] \times [0,1.309017]$, an AABB $\not\subseteq$ RAABB with a volume of less than $5.5451$ cubic units for any sufficiently small $\varepsilon$ (i.e., $\varepsilon \leq 1.1936 \cdot 10^{-7}$ will always work). Of course, this does not prove that $\overline{\overline{P}}_{2,2,2} \left(0^+ \right)$ is an optimal uncrossing covering path for the AABB volume, but constitutes a good example of the variety of hard open problems linked to any $G_{n_1,n_2,\ldots,n_k}$ grid, including the basic case $G_{2,2,2}$.


\section{Conclusion} \label{sec:Conc}

This paper ends the trilogy that started with the introduction of the standard clockwise-\linebreak algorithm \cite{9}, able to perfectly solve (with covering trails of minimal link length) the generalization to $k$-dimensions of Loyd's \textit{nine dots puzzle} \cite{10}.

In the present work, we set a unique length class for every edge of any uncrossing covering path for $G_{n_1,n_2,\ldots,n_k}$, imposing also that the link length must be less than the number of points to be visited, and in addition, every covering path is constrained by a given AABB $\subseteq$ RAABB, so we studied different patterns according to the size of the AABB, with the purpose of reducing, for any $k>1$, the link length of the valid paths.

On the other hand, if we discard the RAABB constraint for our uncrossing covering path, even the minimization of the link length applied to $G_{2,2,\ldots,2}$ turns into a difficult challenge for any $k \geq 3$ (see Acknowledgments and Reference \cite{13}), while a different kind of three-dimensional AABB, such as $\left[-\frac{1}{\sqrt{2}},1+\frac{1}{\sqrt{2}} \right] \times \left[0,1+\frac{1}{\sqrt{2}} \right] \times$\linebreak $\left[0,1+\frac{1}{\sqrt{2}} \right]$ $\not\subseteq$ RAABB (i.e., a box with a volume of only $\frac{\left(1+\sqrt{2} \right)^3}{2}<7.035534$ units$^3$ but which is not entirely contained in the RAABB), let us speculate how $\left(-\frac{1}{\sqrt{2}},\hspace{0.9mm}0,\hspace{0.9mm}1+\frac{1}{\sqrt{2}} \right)$-$(1,\hspace{0.9mm}0,\hspace{0.9mm}0)$-\linebreak$\left(-\frac{1}{\sqrt{2}},\hspace{0.9mm}1+\frac{1}{\sqrt{2}},\hspace{0.9mm}0 \right)$-$\left(\frac{1}{2},\hspace{0.9mm}\frac{1}{2},\hspace{0.9mm}1+\frac{1}{\sqrt{2}} \right)$-$\left(1+\frac{1}{\sqrt{2}},\hspace{0.9mm}1+\frac{1}{\sqrt{2}},\hspace{0.9mm}0 \right)$-$(0,\hspace{0.9mm}0,\hspace{0.9mm}0)$-$\left(1+\frac{1}{\sqrt{2}},\hspace{0.9mm}0,\hspace{0.9mm}1+\frac{1}{\sqrt{2}} \right)$ would solve the aforementioned problem if only we could allow self-intersecting for those $6$ links of prescribed length $1+\sqrt{2}$ units.

Thanks to the M$\Lambda$I-algorithm, we showed that it is possible to provide acceptable covering paths, inside the MAABB, whose edges belong to the same integer (and not only irrational \cite{3}) length class, in two, three, and more dimensions. We hope that the arising questions concerning how many valid length classes there are for any $k$-tuple $(n_1,n_2,\ldots,n_k)$, under many different sets of constraints \cite{2} (even allowing trails, walks or trees \cite{11} instead of merely paths), will provide a challenging environment for new researches in this subfield of graph theory, born from an ancient two-dimensional puzzle which was first published more than one hundred years ago \cite{10}.

Let us now conclude our research on this topic, taking leave of us with a remarkable sentence by a revolutionary mathematician who died in the XIX century:

\noindent \textit{``Quand la concurrence, c'est-\`a-dire l'\'ego\"isme, ne r\`egnera plus dans la science, quand on s'associera pour \'etudier, au lieu d'envoyer aux acad\'emies des paquets cachet\'es, on s'empressera de publier ses moindres observations pour peu qu'elles soient nouvelles et en ajoutant: \linebreak$\ll$ Je ne sais pas le reste.$ \gg$''} \'Evariste Galois, De Ste P\'elagie, December 1831 \cite{12}.


\section*{Acknowledgments}

We sincerely thank Koki Goma from Japan for having discovered and shared, in August 2021, the first known (to the best of our knowledge) unconstrained self-intersecting covering path for $G_{2,2,2}$ with $6$ lines, allowing us to find the related covering paths which have been included in Section \ref{sec:Conc}, and inspiring also the covering circuits represented in Figures \ref{fig:Figure_7c}\&\ref{fig:Figure_8c}.

\makeatletter
\renewcommand{\@biblabel}[1]{[#1]\hfill}
\makeatother

\end{document}